# Enhancing Primary Teacher Training through Academic Portfolios in Advanced Mathematics Courses

Carlos Rojas Bruna

Assistant professor

Pontificia Universidad Católica de Chile

Av. Vicuña Mackenna 4860, Macul, Región Metropolitana

+(56) 2354 4511

ORCID: 0000-0003-0845-9986

carojasb@uc.cl

**Abstract**

The gap between theory and practice in mathematics education, particularly in primary-teacher education, necessitates innovative teaching methodologies.

This paper explores the implementation of academic portfolios as a teaching innovation in Algebra and Number Systems I and II courses within the primary teacher education programme at Pontificia Universidad Católica de Chile.

The methodology involved integrating academic portfolios to align course content with essential learning outcomes for future teaching roles. Implementation begins with a negotiation between students and teachers to establish a learning contract, followed by an overview of course rules, content, objectives, materials, and grading rubrics.

Preliminary findings indicate that this innovative method enhances engagement with mathematical concepts, improves assessment efficacy in teacher training, and may contribute to enhanced preparation of primary mathematics teachers. The study highlights the role of portfolios in making students active participants in their learning, significantly enhancing the educational experience of teacher candidates.

These findings suggest a promising avenue for future educational assessments and methodologies in mathematics, indicating that academic portfolios can bridge the gap between theoretical knowledge and practical applications in mathematics teacher education, potentially enhancing teacher preparation. While this study shows promising results, further research with larger samples and longer timeframes would be beneficial to establish causality and long-term impacts.

Keywords: Mathematics Education Assessment, Teacher Training Outcomes, Academic Portfolio, Teaching Innovation

**Introduction**

The gap between traditional teacher training methods and innovative practices is a critical issue in education. Traditional training often focuses on established practices that may not sufficiently prepare educators for a rapidly changing educational landscape. Throughout this paper, we refer to university students training to become teachers as 'teacher candidates' or 'preservice teachers' to distinguish them from the primary school students they will eventually teach.

By contrast, innovative practices in teacher training emphasise the adoption of new pedagogical methods, technologies, and approaches to address diverse learning needs and enhance student engagement and outcomes. As shown in Figure 1, bridging this gap requires integrating innovative assessment approaches with traditional content knowledge.

[Figure 1]

Numerous studies underscore the significance of bridging this gap. For example, Sailors & Hoffman (2019) discussed how hybrid spaces can assist novice teachers in integrating innovative practices into their teaching by merging academic coursework with practical experiences. Similarly, Stefanova et al. (2019) emphasised the advantages of participatory models in identifying and assessing teacher competencies for open and enquiry-based learning, which can enhance readiness for innovative practices in the classroom.

Innovative teacher-training approaches, such as the STEAM approach discussed by Budarina et al. (2022), concentrate on incorporating practice-oriented methods and modern educational environments to prepare teachers for innovative teaching practices. Süer & Oral (2021) discovered that classroom teachers frequently employ innovative pedagogical practices to align with learner-centred educational approaches, underscoring the importance of integrating innovation into teaching.

This academic portfolio represents a pedagogical innovation with transformative potential in the field of mathematics teacher education. It diverges from traditional evaluative approaches by offering a dynamic compilation of evidence that showcases an individual's learning journey including reflections, achievements, and areas of growth. As a multifaceted instrument, the portfolio not only captures the richness of preservice teachers' academic experiences but also facilitates a nuanced assessment of their progress. The concrete aggregation of educational artefacts allows teacher candidates to engage deeply with mathematical concepts, reflect on their own learning, and develop their professional skills continuously. Thus, the academic portfolio has the potential to become a cornerstone in reimagining how we approach the training and development of aspiring mathematics educators with the aim of cultivating a more reflective, competent, and responsive teaching force.

The central challenge at the heart of the Algebra and Number Systems I and II courses at the Pontificia Universidad Católica de Chile is the evident disconnection between the theoretical rigor of mathematical concepts and their practical applicability in elementary classrooms. This disjunction, which is particularly evident in abstract domains such as number theory, logic, and axiomatic set theory, has necessitated a pedagogical pivot towards methodologies that can bridge this divide. Academic portfolios represent a compendium of knowledge and experience designed to serve as a conduit between theory and practice. This approach was designed to anchor abstract mathematical formalism to the tangible realities of primary mathematics education. The portfolio approach is designed to encourage teacher candidates to actively engage with the material, fostering a deeper conceptual understanding that directly maps onto the skills and insights necessary for their future professional praxis.

The methodology employed experimental design centered on the implementation of an innovative teaching model using academic portfolios. This approach was based on a foundational learning contract established between teacher candidates and instructors. The contract outlines the key components of the portfolio, ensuring mutual accountability and a vested interest in the educational process.

The portfolio serves as a dynamic document capturing iterative reflection and evidence of understanding, providing both formative and summative insights into teacher candidate development. The model requires careful consideration of course rules, a deep dive into class content, and thorough orientation of materials. Explicit grading rubrics articulate the criteria for success, demystify the assessment process, and clarify expectations.

Through this model, preservice teachers engage in the dual role of being learners of advanced mathematics and future educators. This experience was designed to resonate with the impending professional responsibility.

It is imperative that teacher formation in mathematics be research-driven, particularly in a domain in which empirical data and reflective practice combine to shape an effective pedagogy. When academic portfolios are embedded within the curriculum design of teacher training, we transcend the traditional pedagogical constraints. This shift was predicated on the alignment of theoretical frameworks, teaching methodologies, and assessment strategies to create a robust learning environment. Specific examples of these alignments include constructivist learning approaches that emphasize active knowledge construction, formative assessment strategies that provide ongoing feedback, and reflective practice methodologies that encourage metacognitive development. In this context, the portfolio is not merely a repository of completed tasks; it is a dynamic tool that fosters reflective practitioners, who can critically examine the nuances of mathematical instruction. Thus, the integration of academic portfolios is a deliberate design choice informed by research insights that underscore the need for teachers to be adaptable, innovative, and capable of bridging the gap between mathematical theory and pedagogical application. This forward-thinking strategy ensures that future educators are adept at disseminating knowledge and stimulating primary students' intellectual curiosity. Such an approach revitalises the teacher-training landscape, ensuring that it remains at the cutting-edge of educational trends and research imperatives.

Assessment methods in mathematics education have traditionally focused on evaluating students' ability to recall and apply formulas through standardised tests, quizzes, and exams. These traditional assessments often emphasise rote memorisation and procedural knowledge with a primary focus on obtaining correct answers. In this traditional model, teacher candidates are typically passive recipients of knowledge, and assessments are usually summative, occurring at the end of a learning period to measure outcomes. However, in recent years, there has been a significant shift towards innovative assessment methods that aim to engage preservice teachers more actively in the learning process. These innovative assessments focus not only on candidates' knowledge, but also on their critical thinking, problem-solving skills, and creativity. Table 1 compares the key aspects of the traditional and innovative assessment methods in mathematics education, highlighting the stark contrast between these approaches and their potential impact on teacher preparation.

[Table 1]

Despite the growing body of research on portfolio assessment in various disciplines, there exists a conspicuous gap in the context of portfolio utilization in mathematics education. Although portfolios have been extensively documented as reflective assessment tools in a myriad of disciplines, their specific application in mathematics courses, particularly for teacher preparation, remains underexplored. Scholarly investigations into portfolio use have traditionally gravitated towards the humanities and arts, where the interpretation of subjective content aligns with the portfolio's reflective nature (Crowley & Dunn, 1995). However, in the domain of mathematics, the binary nature of mathematical correctness appears at odds with the subjective and holistic evaluations facilitated by portfolios. This dichotomous tension may have contributed to the scarcity of comprehensive studies examining the integration of portfolios as purposeful and evaluative tools in the training of primary school mathematics teachers. It is this lacuna in the research that our paper aims to address.

The preliminary results of integrating academic portfolios as pedagogical methodology in the training of primary school mathematics teachers are encouraging. These results warrant further investigation. Over the two years of implementing this innovation, a pattern of increased teacher candidate engagement and intrinsic motivation became apparent.

Participating preservice teachers exhibited deeper immersion in the subject matter, as evidenced by their sustained efforts and proactive contributions to class discussion. Notably, the sentiments expressed in the reflective pieces within the portfolios indicate a growing sense of ownership of the learning process.

Additionally, anecdotal evidence from instructors suggests improved academic performance and a decrease in absenteeism, indicating a noticeable shift in preservice teachers' attitudes towards advanced mathematics courses. These initial observations have been instrumental in formulating the main thrust of this study, which will proceed to a deeper exploration of the sustained impact of academic portfolios on teacher candidate motivation, understanding of mathematical concepts, and future pedagogical practices.

This study is primarily concerned with anticipating tangible educational benefits that may arise from the integration of an academic portfolio within the teacher training curriculum. Among the expected outcomes is a notable surge in teacher candidate motivation and engagement. This response stems from the portfolio's involvement in learning as an active, reflective process rather than a passive, evaluative one. The assignment of the role of architects of their own knowledge to preservice teachers, coupled with the documentation of their intellectual journeys and pedagogical discoveries, is hypothesised to instil a deeper connection with the subject matter. This heightened engagement is projected not only to enrich the immediate learning experience but also to potentially influence their future classroom practices. As novice educators collate, arrange, and reflect on their portfolios, they simultaneously develop innovative teaching practices. These practices may be transferred to their own teaching, thereby potentially influencing the next generation of mathematicians.

The significance of this study extends beyond the boundaries of mathematics teacher training and encompasses a broader narrative of educational reform. The quality of teacher formation represents a pivotal stage in the educational process, characterised by inherent challenges that can impede future educators' efficacy in the classroom. By investigating the implementation of academic portfolios in the context of primary mathematics teacher training, this research aims to contribute valuable insights into mitigating long-standing issues that surface during the teacher formation stage. These issues include rote memorisation of content that is disconnected from practical teaching scenarios, insufficient development of critical reflective practices, and a lack of engagement and motivation among prospective teachers. By reimagining assessment through a purpose-driven lens, this study postulates that academic portfolios may provide an authentic form for evaluating teacher candidates' abilities to integrate theory with practice, foster deeper commitment to the profession, and lay a stronger foundation for lifelong learning.

Traditional assessment methods, which are commonly employed in mathematics education for teacher training, have frequently been the subject of criticism for their emphasis on rote memorisation and the perceived disconnection between theoretical understanding and practical application. In contrast, academic portfolios represent an innovative paradigm that emphasises the reflection, integration, and application of knowledge. This enables preservice teachers to actively engage with and personalise their learning experiences, potentially reducing the perceived gaps between the content learned in the classroom and its application in a teaching context. The portfolio extends beyond a mere collection of completed tasks; it serves as a curated representation of a teacher candidate's evolving pedagogical identity, showcasing their ability to translate abstract mathematical concepts into tangible, accessible learning experiences for future primary students. Through this holistic approach, the academic portfolio may illuminate the interconnectedness of mathematical ideas and their relevance to real-world teaching scenarios, potentially facilitating the transition from learner to educator.

In conclusion, the academic portfolio is not merely an assessment mechanism; it represents a potentially valuable instrument for enhancing the pedagogical landscape in primary-school mathematics education. Its strategic incorporation into teacher training programs offers an alternative to conventional metrics of evaluation, fostering an environment in which future educators can construct a reflective and analytical narrative of their learning journey. The potential value of the portfolio lies in its dual role as a repository of academic accomplishment and as a reflective prism through which teacher candidates can envision and prepare for the multifaceted challenges of classroom instruction. By exploring this innovative tool, we investigate an approach to the preparation of mathematics teachers that may connect measures of success to authentic teaching competencies and the cultivation of lifelong passion for mathematics education.

**Literature review**

Traditional assessments in teaching mathematics typically involve methods such as standardised tests, quizzes, and exams, which focus on evaluating students' ability to solve mathematical problems and apply formulas (Sarwar et al., 2012). These assessments often emphasise rote memorisation and procedural knowledge, with a primary focus on correct answers rather than problem-solving (Lessani et al., 2017). In traditional assessment methods, teacher candidates are usually passive recipients of knowledge, and the assessment is often summative, occurring at the end of a learning period to measure outcomes (Sarwar et al., 2012).

On the other hand, innovative assessment methods in mathematics teaching aim to engage teacher candidates actively in the learning process and assess not only their knowledge but also their critical thinking, problem-solving skills, and creativity (Feder & Cramer, 2023). These methods include project-based assessments, real-world applications, performance tasks, and portfolios that allow preservice teachers to demonstrate a deeper understanding of mathematical concepts. According to Feder and Cramer's (2023) systematic review of 246 studies, portfolios represent a promising approach that can serve various functions in teacher education, from promoting reflection to facilitating assessment. Innovative assessments focus on formative strategies that provide ongoing feedback to teacher candidates to support their learning processes.

Teachers who adopt innovative assessment strategies in mathematics education often encourage teacher candidates to explore multiple methods for solving problems, engage in discussions, and provide explanations for their reasoning (Gabriele & Joram, 2007). These methods align with constructivist principles, in which preservice teachers are actively involved in constructing their understanding of mathematical concepts (Anderson & Piazza, 1996). Innovative assessment methods can also integrate technology, physical activity, and real-life applications to make mathematics more engaging and relevant to future teachers (Hraste et al. 2018; Yiting et al. 2022).

The quest for excellence in mathematics education places significant emphasis on the professional development of teachers, who are fundamental to shaping the mathematical competence of future generations. Teacher training programs serve as crucibles in which pedagogical skills are honed and theoretical knowledge is fused with practical application. Within this landscape, assessment methods are of paramount importance, influencing not only the acquisition of content knowledge, but also the development of teaching practices that can engage and inspire students. In recent discourse, academic portfolios have emerged as a novel assessment tool that provides a comprehensive view of teacher candidates' abilities and progress (Feder & Cramer, 2023). According to their systematic review of 246 studies, portfolios can serve various functions in teacher education, from promoting reflection to facilitating assessment. Academic portfolios stand at the juncture of reflective practice, methodical assessment, and personalised learning, elements that are increasingly considered vital to effective teacher education. The literature review explores the enactment and efficacy of the academic portfolio as an assessment paradigm, offering insights into its role in contemporary mathematics teacher training, and how it may address pressing needs within the domain.

The theoretical framework of this study rests on the cornerstone concepts of self-regulated learning, motivation, and metacognition, which are critically linked to the effectiveness of mathematics teacher education. Boekaerts (1997) articulated the importance of self-regulated learning as an educational ideal that encompasses a teacher candidate's ability to plan, monitor, and assess their own learning process. This self-directed approach to learning is especially salient in the context of primary school mathematics education, where preservice teachers are required to not only comprehend complex mathematical concepts, but also imbue their future students with the same capacity for self-regulation. Motivation, as explored by Perry, Phillips, and Dowler (2004), serves as a dynamic engine driving teachers' engagement with mathematical content, an essential element for the development of sustained learning and teaching strategies. Finally, Dignath et al.'s meta-analysis (Dignath, Buttner, & Langfeldt, 2008) illuminates effective practices for fostering self-regulated learning strategies among primary school students, highlighting the need for teachers who are adept at teaching these strategies. Such an integrated theoretical schema underscores our endeavour to explore the academic portfolio as an assessment tool fashioned to refine primary mathematics teachers' training by enhancing self-regulation, motivation, and metacognitive awareness.

The literature on perceptions of mathematics education unveils a plethora of psychological elements that influence both teacher and student dynamics within the learning process. Ashcraft and Kirk (2001) contribute to this discourse by illustrating the complexities of math anxiety, a psychological factor that significantly hinders performance. Their work delved into the interplay between working memory and anxiety, effectively positioning math anxiety as a critical barrier to mathematics engagement. This emotional response not only affects students, but can also extend to preservice teachers, potentially perpetuating a cycle of anxiety and underperformance in future generations. On the other hand, Pieronkiewicz (2015) takes a holistic view of the affective domain, acknowledging the multifarious emotional experiences that underpin the process of learning mathematics. Referred to as "affective transgression", this concept encapsulates the challenge of transforming perceived negative emotions into positive and proactive attitudes towards learning. Pieronkiewicz's research underscores the importance of addressing these affective dimensions in teacher training programs to foster environments conducive to mathematical exploration and discovery. These insights highlight the significance of integrating psychological considerations into the fabric of mathematics education, offering a compelling rationale for the need to tailor assessment methods and instructional strategies to ameliorate emotional impediments. Academic portfolios may offer a potential pathway to address these affective concerns by providing preservice teachers with opportunities to reflect on and process their emotional experiences with mathematics in a structured yet personalized format.

The current pedagogical landscape in mathematics teacher training is fraught with myriad assessment tools, each purporting to measure disparate facets of pedagogical competencies and content mastery. Within this diverse array of instruments, academic portfolios have gained prominence, owing to their unique configuration and potential impact on teacher development. Boekaerts' (1997) seminal work on self-regulated learning vividly illustrates portfolios' capacity to foster reflective practitioners who continuously appraise and guide their learning processes. Unlike traditional assessment methods, such as standardised testing, which may provide only a snapshot of a teacher candidate's proficiency, portfolios offer a continuous narrative of growth, encapsulating the evolution of a teacher's thinking, understanding, and instructional planning. Mary L. Crowley & Ken Dunn (1995) accentuate the portfolio's role as a comprehensive record that facilitates a longitudinal perspective on teacher trainees' learning and development. This longitudinal approach is not only distinctively suited for illuminating the depth and breadth of a candidate's mathematical acumen but also for capturing the nuanced interplay of pedagogical skills and content knowledge intricately woven throughout their training experience. Feder and Cramer (2023) further substantiate this view in their systematic review, noting that portfolios demand active involvement from learners, ensuring that teacher candidates are not passive recipients in the evaluative process but rather critical contributors to the discourse of their educational journeys. This evaluative paradigm, therefore, aligns with the modern educational demand for assessments that not only measure outcomes but also engage learners in meaningful, self-regulated learning experiences.

When considering the academic portfolio as a purpose-driven assessment tool, it is essential to explore its efficacy in nurturing motivation and deepening the sense of meaning in advanced mathematics courses. Crowley and Dunn's (1995) pioneering work on mathematics portfolios underscores their multifaceted role not only as a summative showcase of teacher candidate achievement but also as an ongoing reflective process that fosters connections between abstract mathematical theories and real-world applications. The portfolio allows for structured yet flexible documentation of individual learning trajectories, enabling preservice teachers to build and evidence their mathematical journeys. Aligning with Kramarski and Revach's (2009) findings, the process of compiling a portfolio invests teacher candidates with a degree of autonomy and agency in their learning experiences, potentially enhancing their intrinsic motivation. Feder and Cramer (2023) note in their systematic review that while perceptions of portfolios are well-documented, more research is needed to explore the specific circumstances under which portfolios can contribute effectively to teacher education. Careful consideration of such assessment tools is thus crucial to ascertain their impact on motivational dynamics in advanced mathematics courses. Challenging the traditional views of assessment, the portfolio prompts a reassessment of how educators gauge and sustain teacher candidate engagement and whether such tools indeed contribute positively to learners' perceptions of the meaning and purpose in their mathematical studies.

**Background**

The academic setting of this study was carefully integrated into the Primary Teacher Education program curriculum at the Pontificia Universidad Católica de Chile. This robust curriculum comprises a comprehensive array of courses that encompass foundational, pedagogical, and disciplinary knowledge, coupled with didactic instruction. Integral to this educational framework is the progressive development of five discrete internships commencing from the first year, thereby facilitating the rich osmosis of theoretical learning into practical educational arenas.

As students advance to the ninth semester of the program, they are presented with the opportunity to delve deeper into areas of specialisation, namely Natural Sciences, History, Geography and Social Sciences, Language, or Mathematics, each with a unique blend of academic rigor and contextual relevance. The archetype of a graduate from this program is distinguished by adept leadership skills and multifaceted competencies, poised to succeed in various professional contexts, including public and private educational centres, educational management, team leadership, as well as research centres and institutions focused on educational material design and policy analysis. In fine, graduates are envisaged to emerge as reflective, creative, and respectful leaders, exuding a proactive service orientation, a propensity for teamwork, and embodying positive leadership attributes.

Within the mathematics specialization track, the focal point of our teaching innovation project is the implementation of academic portfolios in two cornerstone courses: Algebra and Number Systems I (MAT2920) and Algebra and Number Systems II (MAT2925). These courses are integral components of the ninth and tenth semesters, respectively. MAT2920 is structured to provide a rigorous exploration of the mathematical underpinnings of classical algebra, with particular emphasis on different forms of representation, typical misconceptions, and methods of problem formulation and resolution. The objectives were to scaffold the students' understanding of algebra as a unifying thread in mathematics, foster adeptness in verifying mathematical properties in various number systems and lay the groundwork for the application of algebraic reasoning. Moreover, the course aimed to draw parallels between advanced disciplinary knowledge and scholastic mathematical content. Following this foundation, MAT2925 continues the trajectory by deepening students' conceptual grasp of number systems, underscoring the reflective and application-based facets of arithmetical and proportional thought and their implications for well-rounded mathematics education. Collectively, these courses are designed not merely as academic pursuits but as pivotal stages in the teacher education program, aligning closely with the requisites and practical realities that prospective teachers will encounter in their professional lives.

Table 2 provides a detailed breakdown of the topics and subtopics covered in the Algebra and Number Systems I and II courses.

[Table 2]

This detailed course content outline emphasises the comprehensive nature of the curriculum and its alignment with the objectives of the academic portfolio approach.

A pervasive challenge encountered in the training of primary school mathematics teachers is the alignment of advanced theoretical mathematical concepts with their practical applicability to classroom teaching. Teacher candidates frequently grapple with the perceived incongruence between the depth of disciplinary knowledge imparted in their coursework and the requirements for forthcoming teaching responsibilities. This dissonance manifests in preservice teachers' recurring enquiries about the necessity and relevance of in-depth mathematical understanding, when contrasted with the level they are expected to teach. The palpable gap between the abstract rigor of mathematical theory and the pragmatic demands of pedagogy can result in a sense of disconnection and questioning the value of rigorous disciplinary content within their future professional practice. This tension underscores the need for innovative teaching methodologies that meaningfully bridge theoretical knowledge and practical teacher competencies, thereby addressing the critical need for teacher education programs.

## Methodology/Materials and Methods

### Study Design and Participants

This study employed an experimental design centred on implementing academic portfolios as an innovative teaching methodology in primary mathematics teacher education. The study included a total of 37 teacher candidates enrolled in the Algebra and Number Systems I (MAT2920) and II (MAT2925) courses over four semesters from 2022 to 2023. Specifically, the participants were distributed across four cohorts: 10 students in Algebra and Number Systems I (second semester 2022), 12 students in Algebra and Number Systems II (first semester 2023), 7 students in Algebra and Number Systems I (second semester 2023), and 8 students in Algebra and Number Systems II (first semester 2024). Participants were in their ninth and tenth semesters of the Primary Teacher Education program, specializing in mathematics education. All participants were female, ranging in age from 22 to 25 years, with three semesters of prior classroom experience as teaching assistants (called initial teacher practice) and concurrently completing their professional teaching practice during the study period. Participation in the study was voluntary, and all participants provided informed consent.

**Academic Portfolio Implementation**

The academic portfolio was introduced as a pedagogical innovation aimed at resolving the disconnection between theoretical mathematical concepts and their applications in classroom teaching. Its implementation is congruent with professional development mandates that guide teacher evaluation systems, wherein portfolios constitute key assessment instruments. The portfolio methodology represented [XX%] of the total course grade, with the remaining assessment comprised of [list other assessment components].

The implementation of the academic portfolio commences with a pivotal negotiation phase, in which a foundational learning contract is articulated between instructors and teacher candidates. This contract establishes a bilateral agreement that delineates the explicit expectations, requisites, and guidelines that govern a portfolio's production and assessment. The contract comprehensively outlines the learning objectives, content areas, methodologies for evaluation, deadlines for submissions, and overarching standards. An anonymized example of this learning contract is provided in Appendix B.

Figure 2 illustrates the step-by-step process of portfolio implementation, detailing each phase from initial negotiation to final submission.

[Figure 2]

This flowchart provides a visual representation of the structured approach used to guide teacher candidates through the portfolio creation and submission processes.

**Portfolio Content and Structure**

Algebra and Number Systems I and II courses were meticulously designed to foster a comprehensive understanding of mathematical concepts and their pedagogical applications. Each course session comprised a suite of educational materials including didactic presentations, illustrative examples, and targeted exercises intended to reinforce theoretical understanding and stimulate critical thinking. During these sessions, the preservice teachers were also presented with reflection questions that served to ignite discourse and facilitate the synthesis of disciplinary content with prospective teaching scenarios. Examples of these reflection questions are provided in Appendix C and include prompts such as:

> How might you adapt this algebraic concept for teaching to primary school students?
>
> What connections do you see between this theoretical framework and practical classroom applications?
>
> Identify potential misconceptions students might develop about this topic and how you would address them.

Teaching assistants provided personalized assistance to teacher candidates as they navigated through the course material and the completion of portfolio elements. Teacher candidates submitted initial portfolio drafts, received detailed feedback, and were given opportunities to revise and improve their work before final submission.

The evaluation rubric for academic portfolios was meticulously crafted and served as an explicit guidebook that set forth precise specifications for portfolio assembly and submission. This rubric delineated the format to which the preservice teachers adhered, which included a well-defined structure comprising various critical components. Table 3 provides an overview of the key components of the academic portfolio and criteria used for their assessment.

**[Table 3]**

This table clarifies the specific elements required in a portfolio and the standards against which they are evaluated.

**Data Collection**

Data were collected through multiple sources to provide a comprehensive understanding of the impact of academic portfolios:

- Course Records: Attendance records, digital platform interaction logs (CANVAS), and academic performance data (grades) were collected and analysed to track engagement and achievement patterns.
- Self-Assessment Surveys: At the end of each semester, teacher candidates completed self-assessment surveys regarding their portfolio experience. The survey included Likert-scale questions about time spent on portfolio development, satisfaction with grades, and perceived impact on academic behaviours (self-study, time organization, content understanding, and mental health). The complete survey instrument is included in Appendix A.
- Peer Evaluation Surveys: Teacher candidates also completed peer evaluations, assessing their colleagues' portfolios for dedication, quality, and alignment with grading outcomes. These anonymous evaluations provided additional perspectives on portfolio effectiveness.
- Reflective Components: One-page reflective paragraphs (minimum 200 words) included in each portfolio were analysed for insights into teacher candidates' perceptions and experiences.
- Final Presentations/Discussion Tables: In the first semester of 2023, candidates gave short presentations about their portfolios in both in-person and digital formats. In the subsequent semester, discussion tables were implemented where candidates actively participated in dialogue about questions, they proposed related to their portfolio submissions.

All surveys were voluntary and anonymous, with results aggregated to ensure confidentiality. Measures were taken to ensure teacher candidates felt comfortable providing honest feedback without concern for negative implications by having the surveys administered by a third party not involved in grading.

**Data Analysis**

The analysis employed a mixed-methods approach incorporating both quantitative and qualitative data. Quantitative analysis included descriptive statistics of attendance rates, digital engagement metrics, and grade trends across semesters. Qualitative data from reflective paragraphs and survey responses were thematically analysed to identify patterns in teacher candidates' perceptions and experiences.

For the self-assessment and peer evaluation surveys, responses were categorized and quantified to identify trends in time investment, satisfaction levels, and perceived impacts. The Chilean grading scale ranges from 1.0 to 7.0, with 4.0 representing the minimum passing grade and 7.0 indicating excellence, providing context for interpreting the grade-related findings.

**Results**

This section presents the findings from the implementation of academic portfolios in the primary school teacher training program for mathematics. The results derive from a blend of qualitative and quantitative data sources, including attendance records, digital platform analytics, grade comparisons, and survey responses from both teacher candidates and instructors. Our investigation revealed key findings in three main areas: (1) engagement and participation, (2) academic performance, and (3) perceptions of the portfolio methodology.

On the Chilean grading scale of 1.0 to 7.0, where 4.0 represents the minimum passing grade and 7.0 indicates excellence, our analysis tracked performance across multiple semesters to identify patterns coinciding with portfolio implementation. The findings presented below suggest several positive associations between portfolio usage and various educational outcomes, though we acknowledge that multiple factors may contribute to these observed changes.

### Engagement and Participation

Our quantitative analysis of class attendance and participation shows an encouraging trend towards academic commitment among teacher candidates. Throughout the courses employing the portfolio-based methodology, we observed consistent attendance patterns with presence in classes maintained above 70%. This figure represents solid engagement compared to typical attendance patterns in advanced courses at our institution.

Furthermore, we noted no recorded absences for any of the evaluated activities (portfolio submissions, presentations, and discussion tables), indicating consistent participation in assessment components. While these attendance patterns coincided with the implementation of portfolio methodology, it's important to note that other factors such as cohort characteristics, scheduling, and post-pandemic return to in-person learning may have also contributed to these positive attendance trends.

Digital interaction with course materials was also substantial, with an average of 185 page views on CANVAS per week. While high volume of page views indicates frequent access to materials, this alone cannot be interpreted as deeper engagement without additional data. As one teacher candidate noted in their survey response: "I found myself checking the course materials more frequently than in other courses because I needed to continuously integrate new concepts into my portfolio rather than just studying for exams.

**Academic Performance**

Examining teacher candidate performance in the Algebra and Number Systems I (MAT2920) and Algebra and Number Systems II (MAT2925) courses reveals a notable pattern of grade improvements that coincided with the implementation of the academic portfolio methodology. The portfolios were implemented in MAT2920 in the second semester of 2022, and in MAT2925 in the first semester of 2023.

In the MAT2920 course, the average grade in the second semester of 2022, when portfolios were first introduced, was 5.89 (on a scale of 1.0 to 7.0). This represented an improvement from the previous semester's average score of 4.43. The average grade continued to rise in subsequent semesters: 5.71 in the first semester of 2023 and 6.44 in the second semester of 2023.

Similarly, in the MAT2925 course, the average grade in the first semester of 2023 (when portfolios were implemented) was 5.80, compared with the previous semester's average of 5.52. By the second semester of 2023, the average grade increased further to 6.44.

Figure 4 illustrates the trends in the average grades for MAT2920 and MAT2925 across various semesters.

[Figure 4]

While these improvements in grades coincided with the implementation of portfolios, several other factors may have contributed to these changes, including:

- Different student cohorts with varying abilities and backgrounds
- Post-COVID return to more consistent in-person instruction
- Instructors' growing familiarity with the courses
- Refinements in teaching strategies beyond the portfolio implementation

**Portfolio Assessment Methods**

The implementation of academic portfolios involved diverse assessment approaches that provided teacher candidates with multiple ways to demonstrate their understanding. During the first semester of 2023, preservice teachers engaged in short presentations that demonstrated their conceptual comprehension in both in-person and digital formats. This dual modality accommodated different learning preferences and circumstances while developing various professional skills from public speaking to digital communication.

In the subsequent semester, we incorporated discussion tables as an assessment method, implementing a more interactive and Socratic approach. Teacher candidates actively participated in dialogue around questions they proposed in relation to their portfolio submissions, creating an environment that promoted critical thinking and peer learning. Figure 5 shows preservice teachers' engagement during the final portfolio assessment session.

[Figure 5]

As one instructor noted: "The discussion tables created a completely different dynamic compared to traditional exams. Students were genuinely excited to share their insights and build on each other's ideas, demonstrating a deeper level of understanding than what we typically see in written assessments."

These assessment methods emphasized the portfolio's role not just as a collection of work, but as a catalyst for meaningful discussion and knowledge sharing among future educators.

The importance of these modalities extends beyond mere assessment of understanding; they signify a pedagogical pivot towards learning experiences that echo the collaborative and communicative realities of the teaching profession. The phenomenon illustrated is not simply one of comprehension but of transference, the ability of teacher candidates to convey mathematical concepts in a manner that is pedagogically sound as it is scholarly competent.

The application of an academic portfolio as an instructional method has prompted substantive improvements in bridging the knowledge acquired through coursework with its application in real-world teaching scenarios. Anecdotal evidence and thorough observations reported by instructors indicate that students begin to exhibit a grounded connection between theoretical mathematical concepts and practical teacher-related tasks. Through the deliberate planning of class sessions and the structured design of activities, educators were able to infuse a sense of the relevant purpose into the learning materials, transforming teacher candidates' perceptions of the content from abstract theories to practical tools for classroom implementation. This intentional methical approach may have positively influenced teacher candidates' learning, compelling them to contemplate and execute course content with a vigilant eye towards their future roles as educators. The synthesis of coursework with teaching practice not only advances the effectiveness of the academic portfolio methodology but also suggests a potential impact on teacher candidates' comprehension and integration of mathematical concepts, proving instrumental for their professional development as future primary school teachers.

### Self-Assessment Survey Findings

The implementation of the academic portfolio was further validated through a methodical self-assessment survey completed by the teacher candidates. Across the board, the data indicated a substantial time commitment to portfolio development, with no teacher candidate reporting less than two hours of engagement in any given week. Specifically, 28.6% of the respondents spent between two and four hours, an equal proportion dedicated four to eight hours, while the highest bracket, comprising 42.9%, invested over eight hours leading up to the deadline. This distribution underscores the intense and varied levels of involvement demanded of the portfolio process.

Moreover, the portfolios resonated with the teacher candidates as faithful reflections of their academic endeavours throughout the semester; this sentiment was unanimous among all surveyed participants. Satisfaction with the resultant grades presented a mostly positive picture; 85.7% of the teacher candidates agreed that the grades received were aligned with their efforts and final submissions, while 14.3% (representing a single respondent) expressed dissatisfaction with their grade.

In addressing the influence of portfolios on teacher candidates' academic behaviors, self-assessment highlighted four main aspects. A majority noted impactful enhancements in self-study (85.7%) and time organization (100%), suggesting that portfolios act as catalysts in fostering autonomous learning and efficient management of academic responsibilities. Additionally, over 70 percent felt that portfolios bolstered their understanding of course contents and mental health collectively, indicating that portfolios may offer more than just academic gains; they may serve as a buffer against the academic stresses commonly encountered in rigorous subjects such as mathematics. As one teacher candidate reflected in their portfolio: "The process of creating this portfolio helped me organize my thoughts about mathematics in a way I've never experienced before. I'm starting to see connections between abstract concepts and how I might teach them to children someday." These findings reflect the success of portfolios in meeting the diverse educational and well-being needs of future primary school teachers in our mathematics program.

**Peer Evaluation Results**

The peer evaluation survey results provided valuable additional perspectives on the portfolio process. When estimating time investment of their colleagues, peers reported that four teacher candidates spent between four and eight hours on their portfolios, and three dedicated more than eight hours, consistent with self-reported time commitments. All peer evaluators agreed that the level of dedication was evident in the portfolios they reviewed, and they concurred that the awarded grades aligned with the observable effort and quality.

Teacher candidates identified several common strengths and weaknesses in their peers' portfolio submissions, as summarized in Table 4.

[Table 4]

The strengths identified by peers (n=7) included neatness and clear organization (n=6), effective use of technology (n=4), detailed step-by-step explanations (n=5), and thoughtful reflections (n=3). This feedback suggests that teacher candidates developed the ability to critically evaluate educational materials and recognize elements that contribute to effective mathematical communication.

The weaknesses identified focused primarily on presentation issues rather than conceptual understanding, including incomplete exercises (n=2), unclear statements (n=1), poor image quality (n=2), handwriting legibility (n=1), and limited depth in explanations (n=2). These critiques, while acknowledging areas for improvement, centred on formatting and presentation rather than mathematical understanding, suggesting that preservice teachers' overall conceptual grasp of the material was solid.

As one peer evaluator commented: "I was impressed by how my colleague connected abstract algebra concepts to concrete teaching examples. The organization made it easy to follow their thinking, though some of the handwritten portions were difficult to read."

The collected data suggest progress toward the project's aims—chief among them, the enhancement of teacher candidates' engagement and bridging of the notorious gap between academic content and its tangible implementation in teaching practice. An examination of the results from a pedagogical perspective reveals a pivotal shift: teacher candidates are no longer passive recipients of abstract knowledge but active recipients of a meaningful educational trajectory that resonates with their prospective roles as educators. Portfolio methodology, in its versatile embodiment, appears to be associated with improvements in the comprehension of key mathematical concepts. This is not a trivial outcome, given the inherent abstraction present in courses such as Algebra and Number Systems II. teacher candidates now readily engage with the subjects by integrating definitions and principles into a cohesive pedagogical framework, laying the groundwork for future innovations in their teaching endeavours.

The high level of consensus regarding these benefits, while encouraging, suggests that future research should explore more diverse perspectives, including potential challenges and limitations of the portfolio approach. As one teacher candidate reflected in their portfolio: "The process of creating this portfolio helped me organize my thoughts about mathematics in a way I've never experienced before. I'm starting to see connections between abstract concepts and how I might teach them to children someday". Crucially, the survey data suggests possible connections between portfolio use and improvement in key areas targeted by our objectives, namely enhanced student agency in learning, increased pedagogical coherence, and a fortified connection to practical teaching applications. The convergence of these data points with our pedagogical aims underlines the significant strides made towards reforming mathematics teacher education through innovative assessment strategies.

**Impact on Mathematics Anxiety**

Our analysis also examined teacher candidates' perceptions of mathematics anxiety, a significant psychological barrier to engagement and success in quantitative subjects. The end-of-semester survey revealed that only one respondent (8.3%) felt completely identified with the definition of mathematics anxiety, while seven (58.3%) did not identify with it at all, and the remaining four (33.3%) reported moderate levels of anxiety.

Notably, 11 of the 12 teacher candidates (91.7%) indicated that alternative assessment methods, such as academic portfolios, contributed to mitigating feelings related to mathematics anxiety. When asked to explain how portfolios helped reduce anxiety, several themes emerged from their responses:

- The gradual development process reduced pressure compared to high-stakes exams
- The ability to revise work after feedback allowed for learning from mistakes
- The focus on personal growth rather than comparison with peers created a supportive atmosphere
- The integration of reflection allowed for processing emotional responses to challenging content

As one teacher candidate explained: "In traditional mathematics courses, I always felt judged solely on my ability to solve problems quickly under pressure. With the portfolio, I could take time to work through difficult concepts and show how my understanding evolved. This made the whole experience less stressful and more meaningful."

While satisfaction with the portfolio approach was high, teacher candidates also provided constructive feedback for improvement, including suggestions related to:

- More consistency in the grading process (n=3)
- Additional structure for in-class portfolio activities (n=2)
- More detailed initial instruction on portfolio expectations (n=4)
- Increased availability of teaching assistants for portfolio support (n=2)

**Summary of Findings**

The data gathered through attendance records, academic performance tracking, self-assessments, peer evaluations, and qualitative responses suggest several potential benefits of the academic portfolio approach in mathematics teacher education. During the implementation period, we observed increased class attendance, improved grades, substantial time investment in portfolio development, and positive perceptions regarding the portfolio's impact on learning and anxiety reduction.

Teacher candidates reported that the portfolio approach helped them better connect theoretical mathematical concepts with practical teaching applications, encouraged deeper engagement with course content, and provided opportunities for meaningful reflection on their learning processes. The peer evaluation component additionally fostered critical assessment skills valuable for future educators.

While these findings suggest promising associations between portfolio implementation and various positive outcomes, further research with more rigorous controls would be necessary to establish direct causality. The results nonetheless provide valuable insights into how academic portfolios might serve as a bridge between abstract mathematical knowledge and the practical teaching competencies needed by future primary school mathematics teachers.

**Discussion and conclusion**

**Discussion**

As we begin discussing our findings, it is important to restate the purpose of this study: to examine the role of academic portfolios as a purpose-driven assessment tool in mathematics teacher education programs for primary school teachers. Our research questions focused on whether portfolios could enhance engagement with mathematical concepts, improve assessment efficacy, and potentially contribute to teacher preparation.

**Summary of Key Findings**

Our results indicate several potential benefits associated with the implementation of academic portfolios in mathematics teacher education. Teacher candidates maintained consistent attendance above 70%, demonstrated increased digital engagement with course materials, and showed improved academic performance in both courses where portfolios were implemented. Through self-assessment and peer evaluation surveys, teacher candidates reported that portfolios helped them organize their learning, connect theoretical concepts with practical applications, and reduce mathematics anxiety.

These findings align with previous research on portfolio assessment in teacher education. Feder and Cramer (2023) highlighted in their systematic review of 246 studies that while perceptions of portfolios are well-documented, the specific circumstances under which portfolios contribute effectively to teacher education require further investigation. Our study contributes to this growing body of knowledge by examining portfolio implementation in the specific context of advanced mathematics courses for primary teacher education.

### Implications for Teacher Education

The potential implications of our findings for the training of prospective mathematics teachers merit careful consideration. Traditional assessment approaches in teacher education have often emphasized content knowledge acquisition through standardized testing, which may not adequately prepare candidates for the complexities of classroom teaching. Our findings suggest that a portfolio-based approach may facilitate a more holistic development of teacher competencies.

Teacher candidates' feedback indicates that portfolios may help them make connections between abstract mathematical concepts and classroom applications—a persistent challenge in mathematics teacher education. The structured reflection inherent in portfolio creation appears to support the development of pedagogical thinking alongside mathematical understanding, potentially preparing more well-rounded educators.

However, it is essential to acknowledge that these observed benefits may result from multiple factors beyond the portfolio itself. The comprehensive implementation approach, which included negotiated learning contracts, structured feedback, and diverse presentation formats, likely contributed to the positive outcomes observed. This suggests that portfolios may be most effective when embedded within a thoughtfully designed learning environment that supports their purpose and use.

**Relationship with Existing Literature**

Our findings regarding the potential of portfolios to reduce mathematics anxiety are particularly noteworthy. The link between mathematics anxiety and teacher effectiveness has been well-documented (Ashcraft & Kirk, 2001), with implications for how future teachers might address mathematical concepts in their own classrooms. The reflective and iterative nature of portfolio development may provide a less threatening approach to engaging with challenging mathematical content, as suggested by 91.7% of our participants.

The emphasis on portfolio creation as a process rather than a product aligns with constructivist principles in teacher education (Anderson & Piazza, 1996), where knowledge is actively built rather than passively received. This appears to support Kramarski and Revach's (2009) findings regarding the relationship between portfolio development and learner autonomy, as evidenced by our participants' reports of enhanced self-study habits and time management.

While our observations suggest potential benefits of portfolios in mathematics teacher education, they also echo the caution expressed by Feder and Cramer (2023) regarding the need for more rigorous research designs to establish causal relationships between portfolio use and specific learning outcomes. The largely positive perceptions reported by our participants consistent with broader trends in portfolio research, where affective responses tend to be more thoroughly documented than measurable effects on professional competence.

**Limitations**

Several limitations should be considered when interpreting our findings. First, our study involved a relatively small number of teacher candidates from a single institution, limiting the generalizability of results. Future research would benefit from larger samples across diverse teacher education programs.

Second, the absence of a control group makes it difficult to attribute observed improvements directly to portfolio implementation. Changes in academic performance coincided with portfolio introduction but may have been influenced by other factors, including:

- Different cohort characteristics across semesters

- The transition from pandemic-related restrictions to more normal learning conditions
- Increased instructor familiarity with course content over time
- Concurrent refinements in teaching strategies

Third, our reliance on self-reported data for understanding teacher candidates' experiences with portfolios introduces potential response biases. Participants might have been inclined to report positive experiences to meet perceived expectations, particularly given the high visibility of the portfolio initiative within the program.

Finally, while our study captured perceptions and experiences during the course implementation, we lack data on the long-term impact of portfolio use on teacher candidates' subsequent classroom practice. Longitudinal research following graduates into their teaching careers would provide valuable insights into whether the reported benefits translate into enhanced teaching effectiveness.

### Recommendations for Future Research

Based on our findings and identified limitations, we propose several directions for future research:

1. Conduct comparative studies using matched control groups to better isolate the effects of portfolio implementation from other variables.
2. Develop and validate instruments that can more directly measure the relationship between portfolio use and specific teacher competencies in mathematics education.
3. Design longitudinal studies to track how portfolio experiences in teacher education influence subsequent classroom practices and student outcomes.
4. Investigate the optimal design characteristics of portfolios for mathematics education, including the balance between structure and flexibility, assessment criteria, and integration with other course components.
5. Explore potential differences in portfolio effectiveness across diverse teacher candidate populations and mathematical content areas.

Such research would contribute to a more nuanced understanding of how, when, and for whom portfolio assessment might enhance mathematics teacher preparation.

**Conclusion**

Academic portfolios represent a potentially valuable instrument for enhancing the pedagogical landscape in primary school mathematics education. Their strategic incorporation into teacher training programs offers an alternative to conventional evaluation metrics, fostering an environment in which future educators can construct a reflective and analytical narrative of their learning journey.

Our findings suggest that academic portfolios may help address several persistent challenges in mathematics teacher education, including connecting theory with practice, reducing mathematics anxiety, and promoting self-regulated learning. Teacher candidates' engagement with portfolio creation appears to support the development of both mathematical understanding and pedagogical thinking—essential foundations for effective teaching.

Nevertheless, the limitations of our study underscore the need for caution in attributing causality and for continued research to better understand the specific circumstances under which portfolios most effectively contribute to teacher development. As Feder and Cramer (2023) note, despite the widespread adoption of portfolios in teacher education programs, more rigorous empirical evidence is needed to fully understand their impact.

The potential value of the portfolio approach lies in its ability to position teacher candidates as active architects of their professional development rather than passive recipients of knowledge. By exploring this innovative assessment tool, we contribute to ongoing efforts to develop mathematics teacher preparation approaches that connect measures of success to authentic teaching competencies and foster a lifelong passion for mathematics education.

**Figure 1**

*Bridging the Gap Between Traditional Teacher Training Methods and Innovative Practices*

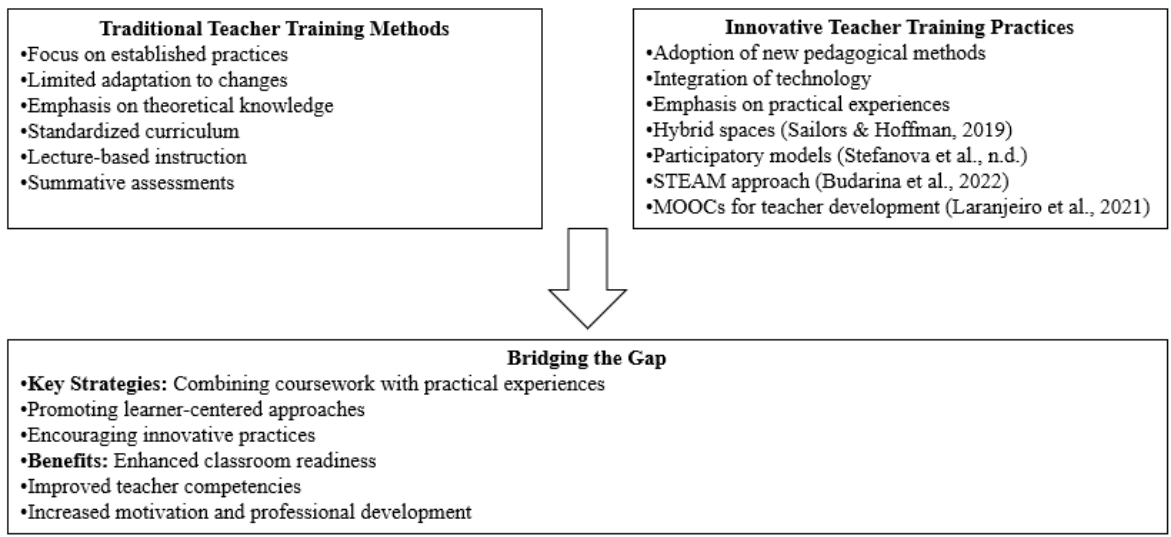

**Table 1**

*Comparison of Traditional and Innovative Assessment Methods in Mathematics Teaching*

| Aspect | Traditional Assessment Methods | Innovative Assessment Methods |
| --- | --- | --- |
| Approach | Standardized tests, quizzes, exams | Project-based assessments, performance tasks, portfolios |
| Focus | Evaluating problem-solving and application of formulas | Assessing critical thinking, problem-solving skills, creativity |
| Student Role | Passive recipients of knowledge | Active participants in the learning process |
| Assessment Timing | Summative (end of learning period) | Formative (ongoing feedback) |
| Emphasis | Rote memorization, procedural knowledge | Deeper understanding of concepts, process of problem-solving |
| Learning Engagement | Limited to correct answers | Encourages multiple problem-solving methods, discussions, explanations |
| Teaching Philosophy | Teacher-centered, focus on correct answers | Student-centred, constructivist principles |
| Integration of Technology | Minimal | Incorporates technology, physical activities, real-life applications |
| Examples of Methods | Standardized tests, quizzes, exams | Real-world applications, MOOC+ flipped classroom, blended teaching |
| Impact on Learning | Often results in superficial understanding | Promotes engagement, critical thinking, and retention of knowledge |

**Table 2**

*Course Content for MAT2920 and MAT2925 - Algebra and Number Systems I and II*

| Course | Topic | Subtopic |
|---|---|---|
| MAT2920 | Algebra and Problem Solving | Variables and expressions; Equations; Inequalities; Algebraic strategies for problem solving |
| | Coordinates, Slopes, and Lines | Rectangular coordinates; Slope; Line equations; Systems of linear equations |
| | Functions and Graphs | Functions; Linear functions; Non-linear functions; Interpretation of graphs |
| | Ratios and Rates | - |
| | Proportions | - |
| | Percentages | Percentages with percentage grids; Calculations with percentages; Approximations |
| MAT2925 | Sets | Basic operations, Cartesian products |
| | Natural Numbers | Peano's axioms; Principle of induction; Operations with natural numbers; Recursive definitions; Order in natural numbers; Principle of well-ordering |
| | Integers | Axioms of integers; Zero; Additive inverses; Operations with signs; Divisibility; Division algorithm; Greatest common divisor; Euclidean algorithm; Least common multiple; Prime numbers; Fundamental theorem of arithmetic |
| | Rational Numbers | Multiplicative inverses; Fractions; Operations with rational numbers; Order in |

| Course | Topic | Subtopic |
|--------|-------|----------|
|        |       | rational numbers; Representation on the number line; Decimal expansion; Operations with decimal expansions |

*Note.* Pontificia Universidad Católica de Chile / January 2016

**Figure 2**

*Overview of the Portfolio Implementation Process*

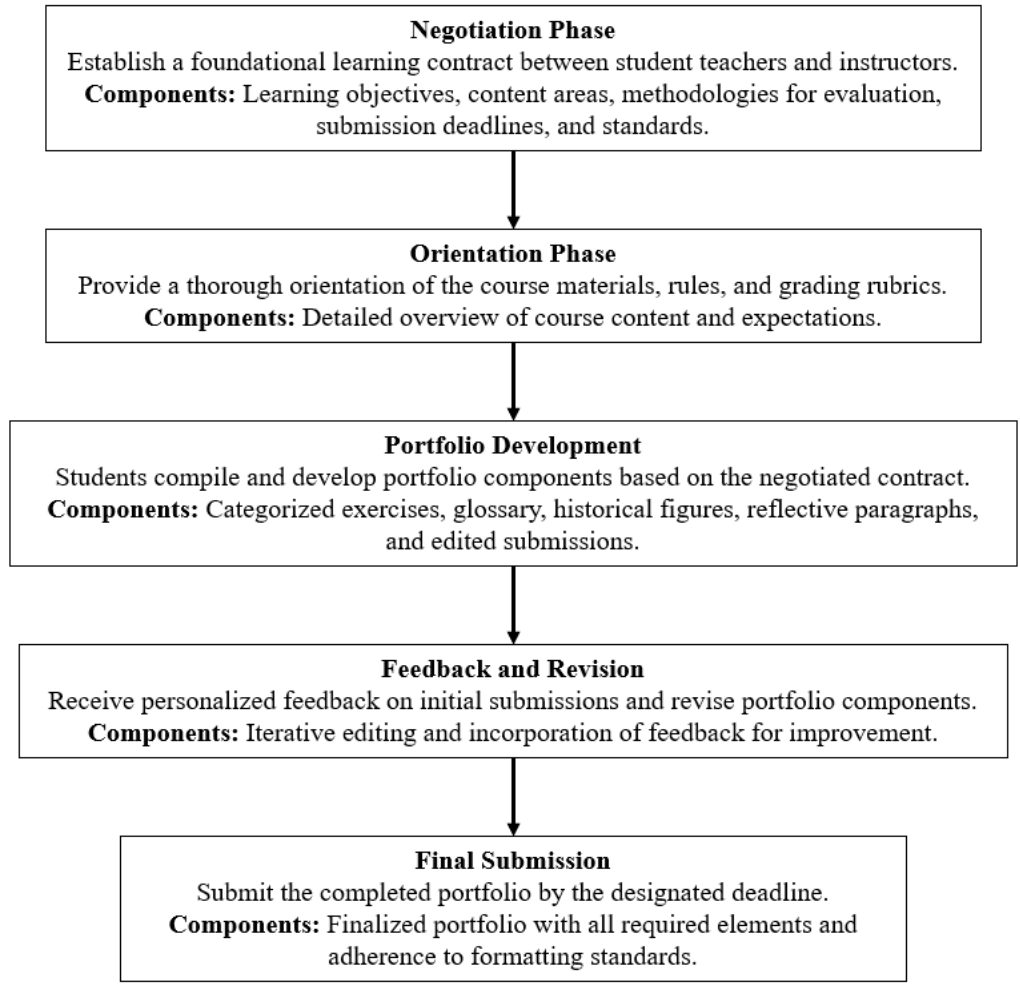

**Table 3**

*Overview of Academic Portfolio Components and Assessment Criteria*

| Component | Description | Assessment Criteria |
|---|---|---|
| Cover Page | Includes student name, instructor's name, course title, and submission date. | Correctness and completeness of information. |
| Categorized Exercises | Responses to exercises categorized and class labeled. | Accuracy, depth of understanding, and completeness of responses. |
| Glossary | Definitions of key concepts encountered during the course, including cited sources. | Accuracy of definitions, proper citation, and relevance to course content. |
| Historical Figures | Chronological list of significant figures in mathematics, providing historical context. | Completeness, relevance, and clarity of the historical context provided. |
| Reflective Paragraph | A minimum 200-word paragraph reflecting on course objectives, content, and its relevance to future teaching practice. | Depth of reflection, relevance to course content, and personal insights. |
| Edited Submissions | Previously submitted assignments edited based on feedback received. | Incorporation of feedback, improvement in clarity and accuracy, and overall enhancement of submissions. |
| Presentation Quality | Adherence to formatting standards, including editing and spelling accuracy. | Presentation quality, adherence to formatting standards, and overall readability. |

| Component | Description | Assessment Criteria |
| --- | --- | --- |
| Timeliness | Submission of portfolio and components by the designated deadlines. | Adherence to deadlines and timeliness of submissions. |

**Figure 3**

*Sample Reflective Paragraph from a Student Portfolio (Original in Spanish)*

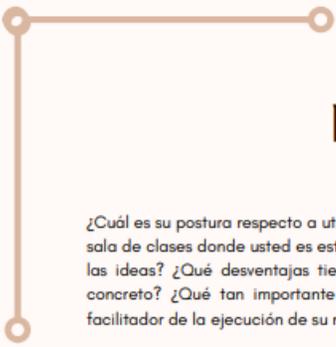

**Figure 4**

*Historical Grades for MAT2920 and MAT2925 Over Semesters*

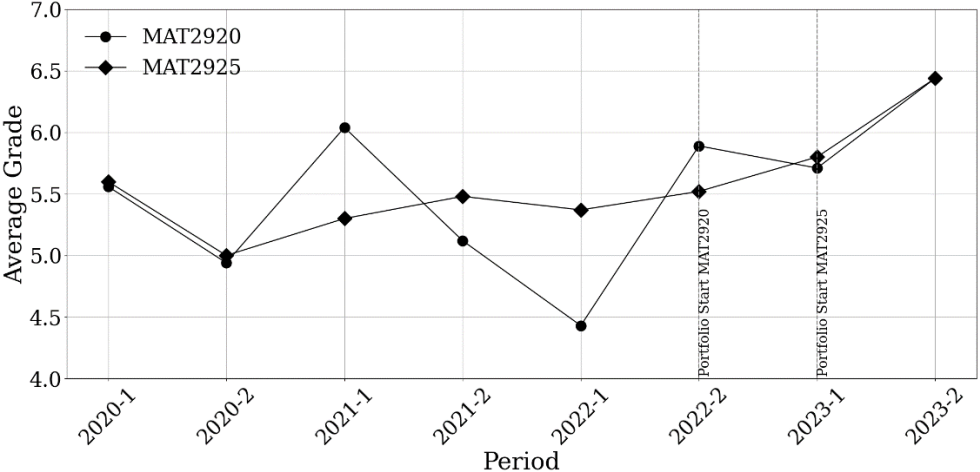

**Figure 5**

*Students Engaging in the Final Exam for Portfolio Submissions*

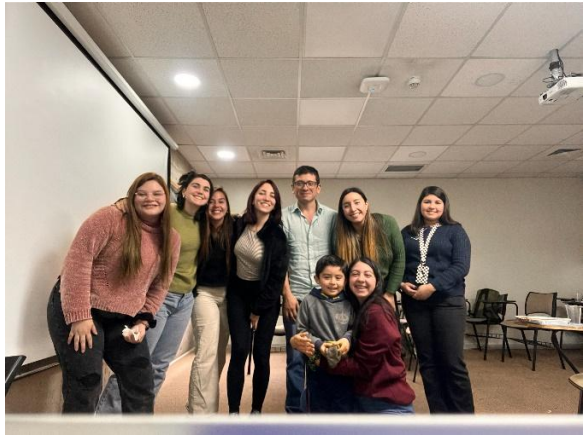

**Table 4**

*Summary of Peer Evaluation Feedback on Portfolio Submissions*

| Category | Strengths Identified by Peers | Weaknesses Identified by Peers |
| --- | --- | --- |
| Content Quality | Detailed explanations, thorough coverage of topics | Incomplete exercises, unclear exercise statements |
| Presentation | Neatness, clear development, and organization | Handwritten legibility, image quality |
| Use of Technology | Effective integration of technology, use of step-by-step explanations | Limited use of digital tools |
| Reflective Elements | Thoughtful reflections, deep personal insights | Superficial reflections, lack of depth in personal insights |
| Timeliness | Timely submission of portfolio components | Occasional late submissions |
| Engagement | Active engagement in the portfolio process, frequent updates | Inconsistent engagement, lack of updates |
| Creativity | Innovative approaches to solving problems, creative presentation of content | Limited creativity, conventional presentation methods |

**Appendix A: Survey Instruments**

**1. Self-Evaluation Questionnaire**

**Spanish (Original)**

**Pregunta 1:** ¿Cuántas horas en promedio dedicó al desarrollo del portafolio durante la semana previa a la entrega?

- de 0 a 2 horas
- de 2 a 4 horas
- de 4 a 8 horas
- más de 8 horas

**Pregunta 2:** Respecto a los temas del curso incluidos para la primera entrega del portafolio. ¿Qué tan seguro se siente respecto al nivel de comprensión alcanzado en estos, pensando en que a futuro podría necesitar enseñarlos en su propia práctica docente? (1 es totalmente inseguro, 5 es totalmente seguro)

- Variables y expresiones
- Ecuaciones
- Desigualdades
- Resolución de problemas

**Pregunta 3:** ¿Cree que su entrega refleja el trabajo realizado durante esta primera parte del semestre?

- Sí
- No
- Parcialmente

**Pregunta 4:** ¿Está de acuerdo con la calificación de su entrega, obtenida a partir de la rúbrica?

- No, debería tener una calificación mayor
- Sí, la calificación es adecuada

**Pregunta 5:** Respecto a la eficiencia en el trabajo realizado. ¿Qué tan eficiente fue su desempeño en cada uno de los siguientes aspectos relacionados con la elaboración de su entrega? (1 es totalmente ineficiente, 5 es totalmente eficiente)

- Organización del tiempo
- Uso de materiales
- Uso de herramientas tecnológicas

**Pregunta 6:** ¿En cuáles de los siguientes factores considera que tiene mayor impacto la utilización del portafolio versus pruebas tradicionales?

- Trabajo autónomo
- Organización del tiempo
- Comprensión de los contenidos
- Salud mental

**English (Translation)**

**Question 1:** How many hours on average did you dedicate to portfolio development during the week prior to submission?

- 0 to 2 hours
- 2 to 4 hours
- 4 to 8 hours
- more than 8 hours

**Question 2:** Regarding the course topics included in the first portfolio submission. How confident do you feel about your level of understanding of these topics, considering you might need to teach them in your future teaching practice? (1 is completely insecure, 5 is completely secure)

- Variables and expressions
- Equations
- Inequalities
- Problem solving

**Question 3:** Do you believe your submission reflects the work completed during this first part of the semester?

- Yes
- No
- Partially

**Question 4:** Do you agree with the grade of your submission, based on the rubric?

- No, I should have received a higher grade
- Yes, the grade is appropriate

**Question 5:** Regarding efficiency in the work performed. How efficient was your performance in each of the following aspects related to the preparation of your submission? (1 is completely inefficient, 5 is completely efficient)

- Time management
- Use of materials
- Use of technological tools

**Question 6:** In which of the following factors do you consider the use of portfolios has a greater impact compared to traditional tests?

- Independent work
- Time management
- Understanding of content
- Mental health

## 2. Peer Evaluation Questionnaire

**Spanish (Original)**

**Pregunta 1:** De acuerdo a lo observado en la entrega de su par, ¿cuántas horas cree que dedicó durante la última semana a la elaboración de la entrega?

- de 0 a 2 horas
- de 2 a 4 horas
- de 4 a 8 horas
- más de 8 horas

**Pregunta 2:** ¿Cree usted que en la entrega de su par se ve reflejada su dedicación para la elaboración del portafolio?

- Sí
- No
- Parcialmente

**Pregunta 3:** ¿Cree usted que la calificación obtenida por su par, es acorde con el trabajo que usted puede apreciar en la entrega?

- Sí, es adecuada
- No, debería tener una calificación mayor
- No, debería tener una calificación menor

**Pregunta 4:** ¿Cuáles son las fortalezas que identificó en el portafolio de su par? (Pregunta abierta)

**Pregunta 5:** ¿Cuáles son las debilidades que identificó en el portafolio de su par? (Pregunta abieerta)

**English (Translation)**

**Question 1:** Based on what you observed in your peer's submission, how many hours do you think they dedicated during the last week to preparing the submission?

- 0 to 2 hours
- 2 to 4 hours
- 4 to 8 hours
- more than 8 hours

**Question 2:** Do you believe that your peer's submission reflects their dedication to the development of the portfolio?

- Yes
- No
- Partially

**Question 3:** Do you believe that the grade obtained by your peer is consistent with the work that you can appreciate in the submission?

- Yes, it is appropriate
- No, they should have received a higher grade
- No, they should have received a lower grade

**Question 4:** What strengths did you identify in your peer's portfolio? (Open question)

**Question 5:** What weaknesses did you identify in your peer's portfolio? (Open question)

## 3. Mathematics Anxiety and Course Evaluation Questionnaire

**Spanish (Original)**

**Pregunta 1:** Lea la siguiente definición. Ansiedad Matemática: "un miedo, tensión y aprensión persistentes relacionados con situaciones que requieren matemáticas (Ashcraft & Kirk, 2001)" ¿Qué tan identificada se siente con lo mencionado en la definición?

- Totalmente identificada
- Parcialmente identificada
- No identificada

**Pregunta 2:** ¿Cree usted que el utilizar instrumentos de evaluación como el portafolio en reemplazo de los tradicionales puede ser un aporte a la mitigación de la aparición de fenómenos como la Ansiedad Matemática?

- Sí
- No

**Pregunta 3:** ¿Qué tan conforme está con respecto al uso del portafolio en el curso como instrumento de evaluación?

- Totalmente conforme
- Parcialmente conforme
- Ni conforme ni disconforme
- Parcialmente disconforme
- Totalmente disconforme

**Pregunta 4:** ¿Qué elementos mejoraría en relación con el desarrollo del curso usando el portafolio como instrumento de evaluación?

- El proceso de calificación
- El desarrollo de las clases
- Las instrucciones
- El desarrollo de las ayudantías

**Pregunta 5:** ¿Le gustaría que el curso Álgebra y Sistemas Numéricos II se trabajara con la misma metodología (mejorando aquellos aspectos perfectibles)?

- Sí
- Me es indiferente
- No

**Pregunta 6:** ¿Está de acuerdo con el que se haga innovación en las metodologías utilizadas en los cursos de su carrera, con el objetivo de mejorar sus aprendizajes y otros elementos?

- Totalmente de acuerdo
- Parcialmente de acuerdo
- Ni de acuerdo ni en desacuerdo
- Parcialmente en desacuerdo
- Totalmente en desacuerdo

**English (Translation)**

**Question 1:** Read the following definition. Mathematics Anxiety: "a persistent fear, tension and apprehension related to situations that require mathematics (Ashcraft & Kirk, 2001)" How identified do you feel with what is mentioned in the definition?

- Completely identified
- Partially identified
- Not identified

**Question 2:** Do you believe that using assessment instruments such as portfolios instead of traditional ones can contribute to mitigating the appearance of phenomena such as Mathematics Anxiety?

- Yes
- No

**Question 3:** How satisfied are you with the use of the portfolio in the course as an assessment instrument?

- Completely satisfied
- Partially satisfied
- Neither satisfied nor dissatisfied
- Partially dissatisfied
- Completely dissatisfied

**Question 4:** What elements would you improve in relation to the development of the course using the portfolio as an assessment instrument?

- The grading process
- The development of classes
- The instructions

- The teaching assistance

**Question 5:** Would you like the Algebra and Number Systems II course to be conducted with the same methodology (improving those aspects that can be perfected)?

- Yes
- I am indifferent
- No

**Question 6:** Do you agree with innovation in the methodologies used in your career courses, with the aim of improving your learning and other elements?

- Strongly agree
- Partially agree
- Neither agree nor disagree
- Partially disagree
- Strongly disagree

**Appendix B: Sample Learning Contract**

**Sample Learning Contract for Algebra and Number Systems Course**

**Spanish (Original)**

**Primera entrega: "Conjuntos y construcción de los números naturales"**
    **Profesor:**
    **Ayudante:**

**¿En qué fecha es la entrega?**

La primera entrega del portafolio es hasta el 17 de abril.

**¿En qué formato debo entregar el portafolio?**

El portafolio puede ser entregado en cualquiera de los siguientes formatos:

- En formato digital. En este caso enviarlo vía correo electrónico.
- En formato no digital. En este caso puede entregarlo al profesor o a la ayudante durante los módulos horarios correspondientes o puede entregarlo durante la semana de entrega en el hall docente de la Facultad de Matemáticas.

**¿Cuáles son los elementos que debe contener el portafolio?**

1. Una hoja inicial donde indique su nombre, el nombre del profesor, el nombre de la ayudante, el nombre del curso, el número de entrega y el título de la entrega
2. El desarrollo de los enunciados etiquetados como "ejercicio" en cada uno de los documentos de las clases, desde la clase 01 a la clase 07 inclusive.
3. Un glosario, al final de los ejercicios, donde enliste el concepto clave en cada una las clases 01 a la 07 con sus respectivas definiciones y la fuente desde donde las obtuvo (sitio web, libro).
4. Una lista con 7 personajes relevantes en relación con los temas desarrollados durante las clases, ordenados cronológicamente e indicando su aporte a la teoría.

5. Un párrafo de reflexión donde usted enuncie y desarrolle una respuesta a la siguiente pregunta: *¿Cuál creo yo que es el propósito de conocer y comprender los fundamentos y definiciones desde el punto de vista formal, de objetos tan elementales como número o conjunto en relación con mi práctica docente? Tomando en consideración que el desarrollo de esta ocurre con un menor nivel de formalidad y abstracción, desde un punto de vista mucho más intuitivo.* El párrafo debe tener una extensión mínima de 200 palabras, puede revisar el sitio web [http://guiastematicas.pucv.cl/subjects/guide.php?subject=c-parrafos-e-acad#tab-0](http://guiastematicas.pucv.cl/subjects/guide.php?subject=c-parrafos-e-acad#tab-0) donde podrá encontrar lineamientos generales acerca de cómo construir el párrafo de manera eficiente.

**Rúbrica de asignación de puntaje**

**Contenido - Hoja inicial (1 punto)**

- Nivel 1 (Excelente - 10 puntos): Están presentes todos los elementos pedidos en la hoja inicial
- Nivel 2 (Muy bueno - 9 puntos): Falta uno de los elementos pedidos en la hoja inicial
- Nivel 3 (Bueno - 8 puntos): Faltan dos de los elementos pedidos en la hoja inicial
- Nivel 4 (Debe mejorar - 6 puntos): Faltan más de dos elementos de los pedidos en la hoja inicial
- Nivel 5 (Malo - 0 puntos): El elemento no está presente en el portafolio

**Contenido - Desarrollo de los ejercicios (4 puntos)**

- Nivel 1 (Excelente - 10 puntos): Están desarrollados de manera correcta los 25 ejercicios
- Nivel 2 (Muy bueno - 9 puntos): Falta, está incompleto o incorrecto el desarrollo de entre 1 y 4 de los 25 ejercicios
- Nivel 3 (Bueno - 8 puntos): Falta, está incompleto o incorrecto el desarrollo de entre 5 y 8 de los 25 ejercicios

- Nivel 4 (Debe mejorar - 6 puntos): Falta, está incompleto o incorrecto el desarrollo de 8 o más ejercicios
- Nivel 5 (Malo - 0 puntos): El elemento no está presente en el portafolio

**Contenido - Glosario y fuente (1 punto)**

- Nivel 1 (Excelente - 10 puntos): Están presentes los 7 conceptos y definiciones claves de la entrega con su respectiva fuente
- Nivel 2 (Muy bueno - 9 puntos): Falta o está incompleta o incorrecta la definición entre 1 o 2 de los conceptos claves de la entrega o no está presente la fuente
- Nivel 3 (Bueno - 8 puntos): Falta o está incompleta o incorrecta la definición de 3 o 4 de los conceptos claves de la entrega o no está presente la fuente
- Nivel 4 (Debe mejorar - 6 puntos): Falta o está incompleta o incorrecta la definición de 5 o más de los conceptos claves de la entrega o no está presente la fuente
- Nivel 5 (Malo - 0 puntos): El elemento no está presente en el portafolio

**Contenido - Lista de personajes relevantes (1 punto)**

- Nivel 1 (Excelente - 10 puntos): Están presentes los 7 personajes ordenados cronológicamente indicando su aporte.
- Nivel 2 (Muy bueno - 9 puntos): Falta o está incompleta o incorrecta la incorporación de entre 1 o 2 personajes en la lista.
- Nivel 3 (Bueno - 8 puntos): Falta o está incompleta o incorrecta la incorporación de entre 3 o 4 personajes en la lista.
- Nivel 4 (Debe mejorar - 6 puntos): Falta o está incompleta o incorrecta la incorporación de más de 5 personajes en la lista.
- Nivel 5 (Malo - 0 puntos): El elemento no está presente en el portafolio

**Contenido - Párrafo libre (2 puntos)**

- Nivel 1 (Excelente - 10 puntos): El párrafo expresa su idea de manera clara y ordenada, dejando, respetando la extensión pedida y las normas de ortografía
- Nivel 2 (Muy bueno - 9 puntos): El párrafo expresa su idea de manera clara y ordenada, sin embargo, no respeta la extensión pedida o algunas normas de ortografía
- Nivel 3 (Bueno - 8 puntos): El párrafo expresa su idea de manera no tan clara y levemente desordenada. Respeta la extensión pedida y las normas de ortografía
- Nivel 4 (Debe mejorar - 6 puntos): El párrafo expresa su idea de manera poco clara y desordenada. No respeta la extensión pedida o incumple algunas normas de ortografía
- Nivel 5 (Malo - 0 puntos): El elemento no está presente en el portafolio

**Presentación - Orden y creatividad (1 punto)**

- Nivel 1 (Excelente - 10 puntos): El portafolio está organizado de manera ordenada. Los elementos están enumerados, se destacan conceptos importantes, las tablas y gráficos son claros y entendibles.
- Nivel 2 (Muy bueno - 9 puntos): El portafolio está organizado de manera ordenada. Pero algunos de los elementos no están enumerados, no se destacan conceptos importantes, Pero las tablas y gráficos son claros y entendibles.
- Nivel 3 (Bueno - 8 puntos): El portafolio está organizado de manera ordenada. Pero algunos de los elementos no están enumerados, no se destacan conceptos importantes y además las tablas y gráficos no son claros ni entendibles.
- Nivel 4 (Debe mejorar - 6 puntos): El portafolio presenta los elementos de forma desordenada. La mayoría de los elementos no están enumerados, no se destacan conceptos importantes y además las tablas y gráficos no son claros ni entendibles.
- Nivel 5 (Malo - 0 puntos): El portafolio presenta los elementos de forma. Los elementos no están enumerados, no se destacan conceptos importantes y además las tablas y gráficos no son claros ni entendibles.

**Presentación - Ortografía (-1 punto)**

- Nivel 1 (Excelente - 0 puntos): El portafolio no presenta errores ortográficos (puntuación, acentuación y gramática)
- Nivel 2 (Muy bueno - 3 puntos): El portafolio tiene menos de 3 errores ortográficos.
- Nivel 3 (Bueno - 6 puntos): El portafolio tiene entre 4 y 6 errores ortográficos.
- Nivel 4 (Debe mejorar - 9 puntos): El portafolio tiene entre 6 y 10 errores ortográficos.
- Nivel 5 (Malo - 10 puntos): El portafolio tiene más de 10 errores ortográficos.

**Penalización - Días de atraso no justificado en la entrega (-1 punto)**

- Nivel 1 (Excelente - 0 puntos): Entregó el portafolio en plazo estipulado.
- Nivel 2 (Muy bueno - 5 punto): Entregó el portafolio con 1 día de atraso.
- Nivel 3 (Bueno - 10 puntos): Entregó el portafolio con 2 días de atraso.
- Nivel 4 (Debe mejorar - 15 puntos): Entregó el portafolio con 3 días de atraso.
- Nivel 5 (Malo - 20 puntos): Entregó el portafolio con más de 3 días de atraso.

**English (Translation)**

**First Submission: "Sets and Construction of Natural Numbers"**
**Professor:**
**Teaching Assistant:**

**When is the submission due?**

The first portfolio submission is due by April 17.

**In what format should I submit the portfolio?**

The portfolio can be submitted in any of the following formats:

- In digital format. In this case, send it via email.

- In non-digital format. In this case, you can deliver it to the professor or the teaching assistant during the corresponding class hours or you can submit it during the delivery week at the faculty hall of the Mathematics Faculty.

**What elements should the portfolio contain?**

1. An initial page indicating your name, the professor's name, the teaching assistant's name, the course name, the submission number, and the title of the submission

2. The development of the statements labeled as "exercise" in each of the class documents, from class 01 to class 07 inclusive.

3. A glossary, at the end of the exercises, listing the key concept in each of classes 01 to 07 with their respective definitions and the source from which they were obtained (website, book).

4. A list of 7 relevant figures related to the topics developed during the classes, arranged chronologically and indicating their contribution to the theory.

5. A reflection paragraph where you state and develop an answer to the following question: *What do I believe is the purpose of knowing and understanding the foundations and definitions from a formal point of view, of such elementary objects as number or set in relation to my teaching practice? Taking into consideration that the development of this practice occurs with a lower level of formality and abstraction, from a much more intuitive point of view.* The paragraph must have a minimum length of 200 words. You can review the website http://guiastematicas.pucv.cl/subjects/guide.php?subject=c-parrafos-e-acad#tab-0 where you can find general guidelines on how to construct the paragraph efficiently.

**Scoring Rubric**

**Content - Initial Page (1 point)**

- Level 1 (Excellent - 10 points): All the requested elements are present on the initial page

- Level 2 (Very Good - 9 points): One of the requested elements is missing from the initial page

- Level 3 (Good - 8 points): Two of the requested elements are missing from the initial page

- Level 4 (Needs Improvement - 6 points): More than two of the requested elements are missing from the initial page

- Level 5 (Poor - 0 points): The element is not present in the portfolio

**Content - Development of Exercises (4 points)**

- Level 1 (Excellent - 10 points): All 25 exercises are correctly developed

- Level 2 (Very Good - 9 points): Between 1 and 4 of the 25 exercises are missing, incomplete, or incorrect

- Level 3 (Good - 8 points): Between 5 and 8 of the 25 exercises are missing, incomplete, or incorrect

- Level 4 (Needs Improvement - 6 points): 8 or more exercises are missing, incomplete, or incorrect

- Level 5 (Poor - 0 points): The element is not present in the portfolio

**Content - Glossary and Source (1 point)**

- Level 1 (Excellent - 10 points): All 7 key concepts and definitions of the submission are present with their respective sources

- Level 2 (Very Good - 9 points): The definition of 1 or 2 of the key concepts is missing, incomplete, or incorrect, or the source is not present

- Level 3 (Good - 8 points): The definition of 3 or 4 of the key concepts is missing, incomplete, or incorrect, or the source is not present

- Level 4 (Needs Improvement - 6 points): The definition of 5 or more of the key concepts is missing, incomplete, or incorrect, or the source is not present

- Level 5 (Poor - 0 points): The element is not present in the portfolio

**Content - List of Relevant Figures (1 point)**

- Level 1 (Excellent - 10 points): All 7 figures are present, arranged chronologically, indicating their contribution
- Level 2 (Very Good - 9 points): The incorporation of 1 or 2 figures in the list is missing, incomplete, or incorrect
- Level 3 (Good - 8 points): The incorporation of 3 or 4 figures in the list is missing, incomplete, or incorrect
- Level 4 (Needs Improvement - 6 points): The incorporation of more than 5 figures in the list is missing, incomplete, or incorrect
- Level 5 (Poor - 0 points): The element is not present in the portfolio

**Content - Reflection Paragraph (2 points)**

- Level 1 (Excellent - 10 points): The paragraph expresses its idea clearly and in an organized manner, respecting the requested length and spelling rules
- Level 2 (Very Good - 9 points): The paragraph expresses its idea clearly and in an organized manner, however, it does not respect the requested length or some spelling rules
- Level 3 (Good - 8 points): The paragraph expresses its idea in a not-so-clear and slightly disorganized manner. It respects the requested length and spelling rules
- Level 4 (Needs Improvement - 6 points): The paragraph expresses its idea in an unclear and disorganized manner. It does not respect the requested length or breaches some spelling rules
- Level 5 (Poor - 0 points): The element is not present in the portfolio

**Presentation - Order and Creativity (1 point)**

- Level 1 (Excellent - 10 points): The portfolio is organized in an orderly manner. The elements are numbered, important concepts are highlighted, and tables and graphics are clear and understandable.

- Level 2 (Very Good - 9 points): The portfolio is organized in an orderly manner. But some of the elements are not numbered, important concepts are not highlighted, but tables and graphics are clear and understandable.

- Level 3 (Good - 8 points): The portfolio is organized in an orderly manner. But some of the elements are not numbered, important concepts are not highlighted, and tables and graphics are not clear or understandable.

- Level 4 (Needs Improvement - 6 points): The portfolio presents the elements in a disorderly manner. Most of the elements are not numbered, important concepts are not highlighted, and tables and graphics are not clear or understandable.

- Level 5 (Poor - 0 points): The portfolio presents the elements in a disorderly form. The elements are not numbered, important concepts are not highlighted, and tables and graphics are not clear or understandable.

**Presentation - Spelling (-1 point)**

- Level 1 (Excellent - 0 points): The portfolio has no spelling errors (punctuation, accentuation, and grammar)

- Level 2 (Very Good - 3 points): The portfolio has fewer than 3 spelling errors

- Level 3 (Good - 6 points): The portfolio has between 4 and 6 spelling errors

- Level 4 (Needs Improvement - 9 points): The portfolio has between 6 and 10 spelling errors

- Level 5 (Poor - 10 points): The portfolio has more than 10 spelling errors

**Penalty - Unjustified Late Days (-1 point)**

- Level 1 (Excellent - 0 points): Submitted the portfolio within the stipulated deadline

- Level 2 (Very Good - 5 points): Submitted the portfolio 1 day late
- Level 3 (Good - 10 points): Submitted the portfolio 2 days late
- Level 4 (Needs Improvement - 15 points): Submitted the portfolio 3 days late
- Level 5 (Poor - 20 points): Submitted the portfolio more than 3 days late

**Appendix C: Sample Reflection Questions**

**Spanish (Original)**

**Algebra y Sistemas Numéricos II**

**1. Escriba un párrafo de reflexión donde usted enuncie y desarrolle una respuesta a la siguiente pregunta:**

"¿Cuál creo yo que es el propósito de conocer y comprender los fundamentos y definiciones desde el punto de vista formal, de objetos tan elementales como número o conjunto en relación con mi práctica docente? Tomando en consideración que el desarrollo de esta ocurre con un menor nivel de formalidad y abstracción, desde un punto de vista mucho más intuitivo."

**2. Escriba un párrafo donde desarrolle la siguiente idea:**

"Al estudiar conceptos matemáticos muchas veces se pierde de vista cual es la definición precisa de los objetos y estos se definen solo mencionando sus aplicaciones, propiedades, características o haciendo alusión a distintas representaciones. Por ejemplo, la definición de la suma de dos números naturales, la definición de que es un número entero, etc. ¿Cree usted que el hecho de al menos conocer las definiciones formales de los objetos es un elemento importante para su práctica docente, al momento de planificar, definir estrategias y presentar de manera apropiada estos conceptos a sus estudiantes? "

**3. Un párrafo libre donde debe desarrollar de la siguiente idea:**

"¿Cuál es su postura respecto a utilización de material concreto en la sala de clases con sus estudiantes y en la sala de clases donde usted es estudiante? ¿Es necesario? ¿Es complementario? ¿Garantiza la comprensión de las ideas? ¿Qué desventajas tiene su uso? ¿Qué sucede con el rol del profesor al trabajar con material concreto? ¿Qué tan importante cree que es el dominio de los conceptos por parte del profesor como facilitador de la ejecución de su rol de conducir la actividad? "

**Algebra y Sistemas Numéricos I**

**4. Un párrafo libre donde usted enuncie y desarrolle su posición respecto a la idea central de esta primera entrega que es:**

"El álgebra es una herramienta para representar información y resolver problemas, ¿es necesaria?, ¿cuáles son las ventajas de su uso?, ¿cuáles son las desventajas?"

**5. Un párrafo libre donde usted enuncie y desarrolle su posición respecto a la idea central de esta primera entrega que es:**

"La representación de problemas de la vida real o matemáticos en el plano cartesiano y la interpretación de la información en el contexto de la situación son una herramienta fundamental para su resolución, así como para la mejor comprensión de los conceptos matemáticos involucrados."

**6. Un párrafo libre donde usted enuncie y desarrolle su posición respecto a lo siguiente:**

"El uso de porcentajes y gráficos para la entrega de información, en diarios, televisión y redes sociales, ¿cree usted que se hace buen uso de ellos? ¿facilitan la comprensión del lector o espectador?"

**English (Translation)**

**Algebra and Number Systems II**

**1. A reflection paragraph where you state and develop an answer to the following question:**

"What do I believe is the purpose of knowing and understanding the foundations and definitions from a formal point of view, of such elementary objects as number or set in relation to my teaching practice? Taking into consideration that the development of this practice occurs with a lower level of formality and abstraction, from a much more intuitive point of view. "

**2. A free paragraph where you must develop the following idea:**

"When studying mathematical concepts, the precise definition of objects is often lost sight of, and these are defined only by mentioning their applications, properties, characteristics, or by alluding to different representations. For example, the definition of the sum of two natural numbers, the definition of what an integer is, etc. Do you believe that at least knowing the formal definitions of objects is an important element for your teaching practice, when planning, defining strategies, and appropriately presenting these concepts to your students?"

**3. A free paragraph where you must develop the following idea:**

"What is your position regarding the use of concrete materials in the classroom with your students and in the classroom where you are a student? Is it necessary? Is it complementary? Does it guarantee understanding of ideas? What disadvantages does its use have? What happens with the role of the teacher when working with concrete materials? How important do you think the teacher's mastery of concepts is as a facilitator of the execution of their role in conducting the activity? "

**Algebra and Number Systems I - Sample Questions**

**4. A free paragraph where you state and develop your position regarding the central idea of this first submission, which is:**

"Algebra is a tool to represent information and solve problems. Is it necessary? What are the advantages of its use? What are the disadvantages?"

**5. A free paragraph where you state and develop your position regarding the central idea of this first submission, which is:**

"The representation of real-life or mathematical problems in the Cartesian plane and the interpretation of information in the context of the situation are a fundamental tool for their resolution, as well as for a better understanding of the mathematical concepts involved."

**6. A free paragraph where you state and develop your position regarding the following:**

"The use of percentages and graphics for the delivery of information, in newspapers, television, and social media. Do you think they are used well? Do they facilitate the reader's or viewer's understanding?"